\documentclass[10pt]{amsart}
\usepackage{amsmath}
\usepackage{amsfonts}
\usepackage{enumerate}
\usepackage{amssymb}
\usepackage{amscd}

\DeclareMathOperator{\gr}{gr}

\DeclareMathOperator{\gld}{gld}
\DeclareMathOperator{\Aut}{Aut}
\DeclareMathOperator{\Hom}{Hom}

\begin{document}
\theoremstyle{plain}
\newtheorem*{trm}{Theorem}
\newtheorem*{lem}{Lemma}
\newtheorem*{prop}{Proposition}
\newtheorem*{thm}{Theorem}
\newtheorem*{thma}{Theorem A}
\newtheorem*{thmb}{Theorem B}
\newtheorem*{thmc}{Theorem C}
\newtheorem*{thmd}{Theorem D}
\newtheorem*{thme}{Theorem E}
\newtheorem*{thmf}{Theorem F}
\newtheorem*{thmg}{Theorem G}
\newtheorem*{thmh}{Theorem H}
\newtheorem*{thmi}{Theorem I}
\newtheorem*{thmj}{Theorem J}
\newtheorem*{example}{Example}
\newtheorem*{cor}{Corollary}
\newtheorem*{conj}{Conjecture}
\newtheorem*{hyp}{Hypothesis}
\newtheorem*{thrm}{Theorem}
\theoremstyle{remark}
\newtheorem{defn}{Definition}
\newtheorem*{rem}{Remark}
\newtheorem*{rems}{Remarks}

\newtheorem*{quest}{Question}
\newtheorem*{notn}{Notation}
\newcommand{\Fp}{\mathbb{F}_p}
\newcommand{\Zp}{\mathbb{Z}_p}
\newcommand{\Qp}{\mathbb{Q}_p}
\newcommand{\Kr}{\mathcal{K}}
\newcommand{\Rees}[1]{\widetilde{#1}}
\newcommand{\invlim}{\lim\limits_{\longleftarrow}}

\title{Primeness, semiprimeness and localisation in Iwasawa algebras}
\author{Konstantin Ardakov and Kenneth A. Brown}
\address{ Department of Mathematics, University of Cambridge;
Department of Mathematics, University of Glasgow, Glasgow G12 8QW}
\email{K.Ardakov@dpmms.cam.ac.uk;
kab@maths.gla.ac.uk}
\begin{abstract}
Necessary and sufficient conditions are given for the completed
group algebras of a compact $p-$adic analytic group with
coefficient ring the p-adic integers or the field of $p$ elements
to be prime, semiprime and a domain. Necessary and sufficient
conditions for the localisation at semiprime ideals related to the
augmentation ideals of closed normal subgroups are found. Some
information is obtained about the Krull and global dimensions of
the localisations. The results extend and complete work of A.
Neumann \cite{Neumann} and J. Coates \emph{et al} \cite{CFKSV}.
\end{abstract}
 \maketitle
\let\le=\leqslant  \let\leq=\leqslant
\let\ge=\geqslant  \let\geq=\geqslant

\section{Introduction} \subsection{} \label{start} In recent years
there has been increasing interest in noncommutative Iwasawa
algebras. These are the completed group algebras
$$ \Lambda_G \quad := \quad \invlim \Zp[G/U], $$
where $\Zp$ denotes the ring of $p-$adic integers, $G$ is a
compact $p-$adic analytic group, and the inverse limit is taken
over the open normal subgroups of $G$. Closely related is the
epimorphic image $\Omega_G$ of $\Lambda_G$,
$$   \Omega_G \quad := \quad \invlim \Fp[G/U], $$
where $\Fp$ is the field of $p$ elements. This increased interest
has been driven by anticipated applications of these rings in
number theory and arithmetic algebraic geometry - see for example
\cite{CFKSV} and the references there for more details.

\subsection{} \label{Iwachat} These algebras were first defined
and studied by Lazard in his celebrated paper \cite{L}. There, he
characterised compact $p-$adic analytic groups as precisely those
topological groups containing an open normal \emph{uniform}
pro$-p$ subgroup of finite index \cite[Corollary 8.34]{DDMS}.
Uniform pro$-p$ groups are defined at \cite[Definition 4.1]{DDMS}.
Lazard showed that when $G$ is uniform $\Lambda_G$ and $\Omega_G$
have the structure of noncommutative formal power series algebras.
In particular, both $\Lambda_G$ and $\Omega_G$ are in this case
\emph{local}, meaning that the factor of each by its Jacobson
radical is simple Artinian, \cite[Exercise 6.2]{DDMS}. (In fact,
$\Lambda_G/J(\Lambda_G) \cong \Omega_G/J(\Omega_G) \cong \Fp.$)
Moreover, (still assuming that $G$ is uniform), $\Lambda_G$ and
$\Omega_G$ are filtered by descending chains of ideals, and the associated graded algebras are
(commutative) polynomial algebras in finitely many variables over
$\Fp$. It follows that $\Lambda_G$ and $\Omega_G$ are Noetherian
domains with good homological properties. A very clear treatment
of all this material can be found in \cite{DDMS}.

If $G$ is now \emph{any} compact $p-$adic analytic group, with
normal uniform open subgroup $U$, say, then clearly
\begin{equation}\label{crossed}
\Lambda_G \quad \cong \quad \Lambda_U \ast (G/U), \end{equation} a
crossed product of the local Noetherian domain $\Lambda_U$ by the
finite group $G/U$. Thus $J(\Lambda_U) \subseteq J(\Lambda_G)$ by
\cite[Theorem 4.2]{Passman}, and hence $\Lambda_G$ is Noetherian
and \emph{semilocal} - that is, $\Lambda_G/J(\Lambda_G)$ is
semisimple Artinian. Similar remarks apply to $\Omega_G$.

\subsection{Primeness and semiprimeness} \label{sprime} Notwithstanding the results
summarised in (\ref{Iwachat}), very little is known about the
(2-sided) ideal structure of noncommutative Iwasawa algebras. The
purpose of this paper is to take the initial steps towards
remedying this. Recall that, if $G$ is any group,
\[\Delta^+(G) = \{x \in G : |G:C_G(x)| < \infty \quad \mbox{and} \quad o(x) < \infty\}.\]
This is a characteristic locally finite subgroup, generated by all
the finite normal subgroups of $G$ \cite[Lemma 5.1]{Passman}. If
$G$ is compact $p-$adic analytic, then $G$ has a torsion-free
subgroup of finite index, since uniform groups are torsion free,
\cite[ Theorem 4.5]{DDMS}, so that $\Delta^+(G)$ is finite, the
unique maximal finite normal subgroup of $G$.

We begin with $\Omega_G$. The first two results are the analogues
for Iwasawa algebras of well-known theorems for ordinary group
algebras, proved by I. G. Connell and D. S. Passman in 1963 and
1962 \cite[Theorems 4.2.10 and 4.2.13]{Pass2}. The proofs for
Iwasawa algebras are however completely different from those for
ordinary group algebras, relying crucially on the crossed product
decomposition (\ref{crossed}).

\begin{thma}
\label{OmegaGPrime} Let $G$ be a compact $p$-adic analytic group. Then
$\Omega_G$ is prime if and only if $\Delta^+(G) = 1$.
\end{thma}

The definition and basic properties of the $p$-\emph{Sylow
subgroups} of a profinite group can be found at \cite[\S 2.2]{W}.

\begin{thmb}
\label{OmegaGSemiprime} Let $G$ be a compact $p$-adic analytic group
and let $P$ be a $p$-Sylow subgroup of $G$. The following are
equivalent:
\begin{enumerate}[{(}i{)}]
\item $\Omega_G$ is semiprime;
\item $\Omega_P$ is semiprime;
\item $\Omega_P$ is prime;
\item $\Delta^+(P) = 1$;
\item $p \nmid |\Delta^+(G)|$.
\end{enumerate}
\end{thmb}

As we show in Theorem F, it is straightforward to deduce the
analogue for $\Lambda_G$ of Theorem A. The analogue for
$\Lambda_G$ of Theorem B, stating that $\Lambda_G$ is always
semiprime for $G$ compact $p-$adic analytic, was noted in
\cite{Neumann} - one simply appeals to the definition of
$\Lambda_G$ as a projective limit of ordinary group rings of
finite groups, over a characteristic zero coefficient ring,
together with the fact that such group rings are semiprime.

\subsection{Zero divisors} \label{zerodiv} We also characterise those $G$ for which $\Omega_G$ is a
domain. The corresponding result for $\Lambda_G$ is due to Neumann
\cite{Neumann}.

\begin{thmc}
\label{OmegaGDomain} Let $G$ be a compact $p$-adic analytic group. Then
$\Omega_G$ is a domain if and only if $G$ is torsion free.
\end{thmc}

\subsection{Localisation at augmentation ideals} \label{loc} A
semiprime ideal $I$ of the Noetherian ring $R$ is said to be
\emph{right localisable} if the subset $\mathcal{C}_R (I) := \{ c
\in R : c + I \textit{ not a zero divisor in } R/I \}$ satisfies
the right Ore condition. Left localisable ideals are defined
analogously, and we say that $I$ is \emph{localisable} if it is
right and left localisable. When the localisation exists, it will
be denoted $R_I$; this is a Noetherian semilocal ring whose
Jacobson radical is $IR_I$. We say that the ideal $I$ satisfies
the \emph{right Artin-Rees property} if, for every right ideal $E$
of $R$, there exists $n \geq 1$ such that $E \cap I^n \subseteq
EI.$ The issue of whether a given ideal is localisable is
intimately connected with the Artin-Rees property - so much so
that, following \cite{Jat}, we shall say that the semiprime ideal
$I$ is \emph{classically right localisable} if it is right
localisable and $IR_I$ satisfies the right Artin-Rees property in
$R_I$; of course ``right" will be omitted where appropriate. It
may well be that the adjective ``classical" is redundant here - in
fact it's a long-standing question of Jacobson whether the
Jacobson radical of a Noetherian ring always has the Artin-Rees
property.

Let $G$ be a compact $p$-adic analytic group and $H$ a closed
normal subgroup of $G$. Since any finite normal subgroup of $H$ is contained in $\Delta^+(H)$ and the product of two finite normal $p'-$subgroups of $H$ is again a finite normal $p'-$subgroup of $H$, we see that $H$ has a largest finite normal $p'-$subgroup, which we will denote by $\Delta^+_{p'}(H)$. This subgroup is clearly characteristic in $H$ and hence is normal in $G$.

Let $w_{H,G}$ denote the kernel of the natural map $\Omega_G \to \Omega_{G/H}$. From our main result on localisation (Theorem D in (\ref{main2})), we can deduce

\begin{thme} Let $G$ be a compact $p-$adic analytic group and let
$H$ be a closed normal subgroup of $G$. Suppose that
$\Omega_{G/H}$ is semiprime. Then the following are equivalent:
\begin{enumerate}[{(}i{)}]
\item $w_{H,G}$ is localisable;
\item $w_{H,G}$ is classically localisable;
\item $H/\Delta^+_{p'}(H)$ is pro$-p$.
\end{enumerate}
\end{thme}

As before, there are versions of Theorems D and E for $\Lambda_G$
also; these are stated and proved in (\ref{lamloc}). Theorem E is
the analogue for compact $p-$adic groups of results of J.E.
Roseblade and P.F. Smith for ordinary group algebras of
polycyclic-by-finite groups; see \cite[Theorem 11.2.15]{Pass2}.

In the case where $H$ is pro$-p$ and $G/H \cong \Zp$, the localisability of the augmentation ideal $v_{H,G} = \ker(\Lambda_G \to \Omega_{G/H})$ of $\Lambda_G$ is proved in \cite[2.4,2.6]{CFKSV}, a result which prompted the work described in this paragraph. Whereas the methods used in that paper are module-theoretic, our methods are more ring-theoretic in nature and are in fact natural extensions of the arguments used by Venjakob in \cite[5.4.1]{Ven}.

There, he proves Theorem E((iii)$\Rightarrow$(i)) and Theorem H((iii)$\Rightarrow$(i)) in the special cases when $H$ is uniform and $G = H \times \Zp$ and when $H=\Zp$ and $G = H \rtimes \Zp$ \cite[Theorem 5.15]{Ven}, and conjectures that these results generalize to the case when $H$ is $p-$valued and $G/H$ is torsionfree in \cite[Conjecture 5.14]{Ven}.

\subsection{Properties of the localisations}

In the last section of this paper, we consider the Krull dimension
and homological properties of the localisations of $\Omega_G$ and
$\Lambda_G$. Let $\Omega_{G,H}$ and $\Lambda_{G,H}$ denote the
localisations of $\Omega_G$ and $\Lambda_G$ at $P_H$ and $I_H$,
respectively. See (\ref{main2}) and (\ref{lamloc}) for the
relevant definitions; we simply note here that if $H$ is a pro-$p$
group with $\Delta^+ (G/H) = 1$ then $P_H = w_{H,G}$ and $I_H
= v_{H,G}.$

\begin{thmi}
Let $G$ be a compact $p$-adic analytic group and let $H$ be a closed normal subgroup of $G$. Then
\begin{enumerate}[{(}i{)}]
\item $\Kr(\Omega_H) \leq \Kr(\Omega_{G,H}) \leq \dim H$,
\item $\Kr(\Lambda_{G,H}) = \Kr(\Omega_{G,H}) + 1$.
\end{enumerate}
\end{thmi}

It is proved in \cite{Ard} that $\lambda(\mathfrak{h}) \leq
\Kr(\Omega_H)$, where $\mathfrak{h}$ is the $\Qp-$Lie algebra of
$H$ and $\lambda(\mathfrak{h})$ is the maximum length $m$ of
chains
$0=\mathfrak{h}_0<\mathfrak{h}_1<\ldots<\mathfrak{h}_m=\mathfrak{h}$
of sub-Lie-algebras of $\mathfrak{h}$. When $\mathfrak{h}$ is solvable,
$\Kr(\Omega_H) =\dim H$ and so (i) consists of equalities, but in
general $\Kr(\Omega_H) < \dim H$.

The parallel result on global dimensions is part of the following
theorem. The unexplained terminology is defined in \ref{hom}. When $H=G$ and $G$ has no elements of order $p$, we obtain an earlier result of Venjakob \cite[Theorem 3.26]{Ven2}.

\begin{thmj}
Let $G$ be a compact $p$-adic analytic group and let $H$ be a
closed normal subgroup of $G$. Then
\begin{enumerate}[{(}i{)}]
\item  $\Omega_{G,H}$ is Auslander-Gorenstein of injective
dimension $\dim H.$ \item $\Lambda_{G,H}$ is Auslander-Gorenstein
of injective dimension $\dim H + 1.$ \item If $G$ contains no
element of order $p$ then $\Omega_{G,H}$ and $\Lambda_{G,H}$ have
finite global dimension. In this case both algebras are
Auslander-regular, and $$\gld(\Omega_{G,H}) = \gld(\Omega_H) =
\dim H,$$ while $$\gld(\Lambda_{G,H}) = \gld(\Lambda_{H}) =
\gld(\Omega_{G,H}) + 1.$$ \item If $\Omega_{G,H}$ or
$\Lambda_{G,H}$ has finite global dimension then the inverse image
in $G$ of $\Delta^+(G/H)$ contains no elements of order $p.$
\end{enumerate}
\end{thmj}

\section{Primeness and semiprimeness of $\Omega_G$}

\subsection{Uniform pro$-p$ groups} \label{uniform}
We begin with some technical results concerning $\Omega_N$ for uniform $N$.
The following notation will be fixed for the rest of this paper:
\[\epsilon = \left\{ \begin{array}{l}2 \quad \mbox{ if } \quad p=2 \\ 1 \quad \mbox{ otherwise.} \end{array} \right. \]

\begin{lem}
Let $N$ be a uniform pro$-p$ group and let $\mathfrak{m}$ denote the maximal ideal of $\Omega_N$. Then
\begin{enumerate}[{(}i{)}]
\item $[\Omega_N, \mathfrak{m}^k] \leq \mathfrak{m}^{k + p^\epsilon - 1}$ for all $k\geq 1$,
\item $(1 + \mathfrak{m}^{p^\epsilon - 1}) \cap N \leq N^{p^\epsilon}.$
\end{enumerate}
\end{lem}
\begin{proof}
To avoid confusion, we will denote the Lie commutator in $\Omega_N$ by $[a,b]:= ab-ba$, and the group commutator in $N$ by $(x,y) := x^{-1}y^{-1}xy$.

Since $N$ is uniform and so in particular powerful, $(N,N) \leq
N^{p^\epsilon}$, \cite[Definition 3.1 and Lemma 3.4]{DDMS}. From
this, it's easy to deduce that $[\Omega_N, \mathfrak{m}] =
[\mathfrak{m},\mathfrak{m}] \leq \mathfrak{m}^{p^\epsilon}$. Part
(i) now follows by an straightforward induction.

Next, $gN^p \mapsto g - 1 + \mathfrak{m}^2$ gives an isomorphism
of the groups $N/N^p$ and $\mathfrak{m}/\mathfrak{m}^2$, so $(1 +
\mathfrak{m}^2) \cap N = N^p$. If $p$ is odd, part (ii) follows,
so assume that $p=2$. Choose a topological generating set $\{a_1,
\ldots a_d\}$ for $N$ and let $g \in (1 + \mathfrak{m}^3) \cap N$.
Since $(1 + \mathfrak{m}^3) \cap N \leq N^2$ by the above, we can
write $g = a_1^{2\mu_1} \cdots a_d^{2\mu_d}$ for some $\mu_i \in
\mathbb{Z}_2$ by \cite[Theorems 3.6(iii) and 3.7 and Proposition
1.28]{DDMS}. Set $b_i = a_i-1 \in \Omega_N$. Then
\[ g - 1 = (1 + b_1^2)^{\mu_1} \cdots (1 + b_d^2)^{\mu_d} - 1
\equiv \bar{\mu}_1 b_1^2 + \ldots + \bar{\mu}_d b_d^2 \quad \mod
\mathfrak{m}^3,\] where $\bar{\mu}_i$ denotes the reduction of
$\mu_i$ modulo 2. By \cite[Theorem 7.24]{DDMS}, the elements
$b_i^2$ are linearly independent modulo $\mathfrak{m}^3$, so
$\mu_i \in 2\mathbb{Z}_2$ for all $i$. Hence $g \in N^4$ as
required.
\end{proof}

Recall that if $N$ is uniform, then $\Omega_N$ is a Noetherian
domain \cite[Corollary 7.25]{DDMS}, and as such has a division
ring of quotients $D_N$ \cite[Theorem 2.1.14]{MCR}. Define the
\emph{degree} $\deg r$ of a nonzero element $r \in \Omega_N$ to be
the greatest non-negative integer $k$ such that $r$ is in the
$k$th power of the maximal ideal of $\Omega_N$. The following
result will be crucial to our proof of Theorem A.

\begin{prop}
\label{OuterAction}
Let $N$ be a uniform pro$-p$ group and let $\varphi \in \Aut(N)$. Suppose that the natural extension of $\varphi$ to $D_N$ is inner. Then $[N,\varphi] \subseteq N^{p^\epsilon}$.
\end{prop}
\begin{proof}
As in Lemma \ref{uniform}, let $\mathfrak{m}$ be the maximal ideal
of $\Omega_N$. By assumption, we can find a nonzero element $x \in
D_N$ such that
\[ x^{-1} g x = g^\varphi \quad \mbox{for all} \quad g \in N.\]
We can also find nonzero elements $a,b,c,d \in \Omega_N$ such that
$x=bd^{-1} = a^{-1}c$; note that this forces \begin{equation}  ab
= cd. \label{vogts}\end{equation} Fix $g\in N$. From the above, it
follows that
\[agb = cg^\varphi d.\]
Now, by the first part of the lemma, and using (\ref{vogts})
twice,
\[
\begin{array}{lll}
c(g^\varphi - g)d &=& cg^\varphi d - (ab-cd)g - cgd \\
                                    &=& a[g,b] - c[g,d] \\
                                    & \in & \mathfrak{m}^{\deg a} \mathfrak{m}^{\deg b + p^\epsilon - 1} + \mathfrak{m}^{\deg c} \mathfrak{m}^{\deg d + p^\epsilon - 1} \\
                                    &=& \mathfrak{m}^{\deg c + \deg d + p^\epsilon - 1}.
\end{array}
\]
Because $\gr_{\mathfrak{m}}\Omega_N$ is a domain \cite[Theorem 7.24]{DDMS}, we see that
\[ g^\varphi - g \in \mathfrak{m}^{p^\epsilon - 1}.\]
Hence by the second part of the lemma, we have
\[ g^{-1}g^\varphi \in (1 + \mathfrak{m}^{p^\epsilon - 1}) \cap N \leq N^{p^\epsilon}\]
for all $g \in N$, so $[N,\varphi] \subseteq N^{p^\epsilon}$ as required.
\end{proof}

\subsection{Primeness - special cases} \label{subpop}
Let $G$ be a compact $p-$adic analytic group, and let $N$ be a closed normal uniform subgroup of $G$. Define
\[E_G(N) = \{x \in G : [N,x] \subseteq N^{p^\epsilon}\}.\]
This is the kernel of the conjugation action of $G$ on the finite
set $N/N^{p^\epsilon}$ and as such is an open normal subgroup of $G$. Since
$N$ is uniform, $N$ is always contained in $E_G(N)$.

We begin the proof of Theorem A by
taking care of two special cases.

\begin{prop}
\label{OuterActionSpecialCase} Let $G$ be a compact $p$-adic analytic
group. Suppose $N$ is an open normal uniform subgroup of $G$ such
that $E_G(N) = N$. Then $\Omega_G$ is prime.
\end{prop}
\begin{proof}
We have observed in the introduction that $\Omega_G$ is a crossed product of the Noetherian domain $\Omega_N$ with the finite group $\overline{G} := G/N$. The set $S = \Omega_N \backslash \{0\}$ is a $\overline{G}$-invariant Ore set in $\Omega_N$ and $D_N=\Omega_NS^{-1}$ is the division ring of fractions of $\Omega_N$. By \cite[Lemma 37.7]{Pass2}, $S$ is an Ore set in $\Omega_N*\overline{G}$ consisting of regular elements, and
\[(\Omega_N*\overline{G})S^{-1} \cong D_N * \overline{G}.\]
Now, if $1 \neq \bar{g} \in \overline{G}$, conjugation by $\bar{g}$ gives an outer automorphism of $D_N$ by Proposition \ref{OuterAction}, because $E_G(N) = N$. By \cite[Theorem 7.8.12]{MCR}, $D_N * \overline{G}$ is a simple ring. Since this is a partial quotient ring of $\Omega_G = \Omega_N * \overline{G}$, it follows that $\Omega_G$ is prime, as required.
\end{proof}

We shall use below the elementary fact that, for any group $G$,
\begin{equation}
 \label{del}
 \textit{  if $H$ is a subgroup of $G$ of finite index, then }
\Delta^+(H) \leq \Delta^+(G). \end{equation}

\begin{lem} Let $G$ be a compact $p-$adic analytic group such that
$\Delta^+ (G)$ is a $p-$group. Let $N$ be a normal open uniform
subgroup of $G$. Then $E_G (N)$ is a pro$-p$ group.
\end{lem}
\begin{proof}
Let $C = C_G(N)$ denote the centralizer of $N$ in $G$, so $C$ is
closed by \cite[Exercises 0.4(2)]{W}. Let $q\neq p$ be a prime and
let $Q$ be a $q$-Sylow subgroup of $C$; since $G$ is virtually
pro$-p$, $Q$ is finite. Since $Q$ centralizes $N$, $Q$ is a finite
normal subgroup of $NQ$ so $Q \leq \Delta^+(NQ) \leq \Delta^+(G)$
by (\ref{del}). Hence $Q = 1$ and $C$ is a pro$-p$ group.

Now, consider the conjugation action of $E = E_G(N)$ on $N$. This
has kernel $C \cap E$ which is a pro$-p$ group since $C$ is, and
the image is contained in
\[\Gamma := \{\varphi \in \Aut(N) : [N,\varphi]\subseteq N^{p^\epsilon}\}.\]
By \cite[Corollary 4.18]{DDMS}, $\Aut(N)$ may be identified with a
subgroup of $GL_d(\Zp)$ where $d = \dim(N)$. It's clear that under
this identification $\Gamma$ is contained in $\{g \in GL_d(\Zp) :
g \equiv 1$ (mod $p^{\epsilon})\}$ which is a pro$-p$ group by \cite[Theorem
5.2]{DDMS}. Hence $\Gamma$ is a pro$-p$ group and thus so is $E$
by \cite[Proposition 1.11(ii)]{DDMS}.
\end{proof}

\begin{cor}
\label{UniformSylow} Let $G$ be a compact $p$-adic analytic group such
that $\Delta^+(G) = 1$. Suppose that a $p-$Sylow subgroup $N$ of
$G$ is normal and uniform. Then $E_G(N) = N$ and $\Omega_G$ is
prime.
\end{cor}
\begin{proof}

Since $E_G(N)$ contains $N$, the lemma shows that $E_G (N) = N$.
Therefore $\Omega_G$ is prime by the proposition.
\end{proof}

\subsection{Proof of Theorem A} \label{thmaproof} The proof needs a technical lemma which will be used to
translate the group-theoretic condition $\Delta^+(G) = 1$ into a
condition involving crossed products. Compare the argument with
the proof of \cite[Lemma 2.2(iii)]{DDMS}.

\begin{lem}
\label{Technical} Let $G$ be a pro$-p$ group of finite rank such
that $\Delta^+(G) = 1$. Suppose $N$ is an open normal uniform
subgroup of $G$ and let $x \in E = E_G(N)$. Then $H = \langle
N,x\rangle$ is uniform.
\end{lem}
\begin{proof}
Since $[H,H] = [N,H] \leq [N,E] \leq N^{p^\epsilon} \leq H^{p^\epsilon}$, $H$ is a
finitely generated powerful pro$-p$ group. By \cite[Theorem 4.20]{DDMS}, the elements of finite
order in $H$ form a characteristic subgroup $T = \Delta^+(H)$ and $H/T$ is uniform. By (\ref{del}), $T = 1$ and the result follows.
\end{proof}

Recall the definition of $w_{H,G}$ (or $w_H$) from
(\ref{loc}).

\begin{proof}[Proof of Theorem A]
Suppose that $\Omega_G$ is prime. Let $e = \sum_{g \in
\Delta^+(G)} g$, a nonzero central element of $\Omega_G$. Since
$w_{\Delta^+(G)} . e\Omega_G = 0$, $w_{\Delta^+(G)} = 0$, and
hence $\Delta^+(G) = 1$.

Now suppose that $\Delta^+(G) = 1$. Choose an open normal uniform
subgroup $N$ of $G$, so that, by (\ref{crossed}), $\Omega_G$ is a crossed product of the Noetherian
domain $\Omega_N$ by the finite group $G/N$. Let $q$ be a prime
(possibly equal to $p$) and let $H/N$ be a $q$-Sylow subgroup of
$G/N$. By \cite[Theorem 17.5]{Passman}, it's sufficient to show
that $\Omega_H = \Omega_N
* (H/N)$ is prime for all such $H$. If $q \neq p$, we see that $N$ is the $p$-Sylow
subgroup of $H$; since $\Delta^+(H)=1$ by (\ref{del}), $\Omega_H$
is prime by Corollary \ref{UniformSylow}.

Suppose therefore that $q = p$, so $H/N$ is a finite $p$-group. By
\cite[Proposition 16.4]{Passman} it's sufficient to prove that
$\Omega_H = \Omega_N * (H/N)$ is semiprime. By \cite[Theorem
18.10]{Passman}, this will follow from the semiprimeness of
$\Omega_P$, where $P$ is an arbitrary subgroup of $H$ such that $N
\subseteq P$ and $P/N$ is elementary abelian.

Let $P$ be any such group, and choose $K$, a maximal uniform
subgroup of $P$ containing $N$. Since $P/N$ is abelian, $K$ is an
open normal uniform subgroup of $P$. If $x \in E_P(K)$, then
$\langle K,x \rangle$ is a uniform subgroup of $P$ containing $K$
by the lemma. The maximality of $K$ forces $x \in K$,
so $E_P(K) = K$. Hence $\Omega_P$ is prime by Proposition
\ref{OuterActionSpecialCase}, as required.
\end{proof}

\subsection{Proof of Theorem B}\label{thmbpf} Choose an
open normal uniform subgroup $N$ of $G$ contained in $P$. Then, as
in (\ref{crossed}), $\Omega_G = \Omega_N * (G/N)$ is a crossed
product of the Noetherian domain $\Omega_N$ by the finite group
$G/N$. It's clear that $P/N$ is a Sylow $p$-subgroup of $G/N$ and
that $\Omega_P = \Omega_N * (P/N)$. The equivalence of (i) and
(ii) now follows from \cite[Corollary 18.11]{Passman}, whereas
(ii) and (iii) are equivalent by \cite[Proposition 16.4]{Passman}.

Now, (iii) and (iv) are equivalent by Theorem A. Since $P \cap
\Delta^+(G)$ is a finite normal subgroup of $P$, $P \cap
\Delta^+(G) \leq \Delta^+(P)$. Since $\Delta^+(P) \leq
\Delta^+(G)$ by (\ref{del}), $P \cap \Delta^+(G) = \Delta^+(P)$.
But $P \cap \Delta^+(G)$ is a $p$-Sylow subgroup of $\Delta^+(G)$
by \cite[Proposition 2.2.3(a)]{W}, and the equivalence of (iv) and
(v) follows.

\subsection{Proof of Theorem C}\label{thmcpf}
Suppose $G$ has torsion. Then we can find $1\neq x \in G$ such that $x^n = 1$. Now $(x-1)(x^{n-1} + \ldots + 1) = 0$ shows that $\Omega_G$ has zero divisors.

Now assume $G$ is torsion free. The proof follows the line of
Neumann's for $\Lambda_G$ \cite{Neumann} in appealing to a result
of Walker \cite{Walker} which guarantees that a semiprime Noetherian
ring $R$ of finite global (homological) dimension, with $R/J(R)$ a
field, is a domain.  For certainly $\Omega_G$ is Noetherian, and
it's semiprime by Theorem B. Any $q$-Sylow subgroup of $G$ with
$q\neq p$ must be finite since $G$ is virtually pro$-p$. Hence $G$
is a pro$-p$ group and therefore $J(\Omega_G) = w_G$. Finally, the fact that $\Omega_G$ has finite global dimension follows from \cite[Theorem 4.1]{Bru} and \cite[Corollaire 1]{Serre}.

\section{Localisable ideals in $\Omega_G$}\label{local}\subsection{}\label {main2}

Let $G$ be a compact $p$-adic analytic group and $H$ a closed normal
subgroup of $G$. Recall the definition of $w_{H,G}$ from
(\ref{loc}). When there is no danger of
confusion, we will shorten $w_{H,G}$ to $w_H$. By \cite[Lemma
4.5]{Bru}, $\Omega_G$ is a flat $\Omega_H$-module, so that
\[ w_{H,G} = w_{H,H}\Omega_G = (H-1)\Omega_G = \Omega_Gw_{H,H}.\]
Note that if $H$ is pro$-p$, then $w_{H,H}$ is the unique maximal ideal of $\Omega_H$.

\begin{lem}
\label{FNPSubgroup}
Let $H$ be a finite normal $p$-subgroup of $G$. Then $w_H$ is nilpotent.
\end{lem}
\begin{proof}
Since $H$ is a finite $p$-group, $\Omega_H = \Fp[H]$ is local
Artinian with maximal ideal $\mathfrak{m} = w_{H,H}$,
 with $\mathfrak{m}^n = 0$ for some $n$, by \cite[Lemma 3.1.6]{Pass2}. Since $H$ is normal,
$w_H^n = (\mathfrak{m}\Omega_G)^n = \mathfrak{m}^n\Omega_G = 0.$
\end{proof}

\subsection{}\label{extra} Now let $H$ be an arbitrary closed normal subgroup of
$G$. By
 \cite[Corollary 8.34]{DDMS}, $H$ contains an open
pro$-p$ subgroup $J$, which can be chosen normal in $G$.

Let $P_H \triangleleft \Omega_G$ be the prime radical $\sqrt{w_J}$
of $w_J$, the smallest semiprime ideal of
$\Omega_G$ containing $w_J$, \cite[0.2.8]{MCR}. By \cite[Theorem 2.3.7]{MCR}, since
$\Omega_G$ is Noetherian, $P_H$ is the biggest ideal of $\Omega_G$
which is nilpotent \textit{modulo} $w_J$. The following lemma is
essentially due to Peter Schneider, \cite[Lemma 2.5]{CFKSV}.

\begin{lem}
\label{PHWellDefined}
$P_H$ is independent of the choice of $J$.
\end{lem}
\begin{proof}
Let $J$ and $K$ be open pro$-p$ subgroups of $H$ which are normal in $G$. Then $J \cap K$ also has these properties,
 so it's sufficient to consider the case when $J \subseteq K$. Now $K/J$ is a finite normal $p$-subgroup of $G/J$,
 so $w_K^n \subseteq w_J$ for some $n$ by Lemma \ref{FNPSubgroup}. Since $J \subseteq K$, $w_J \subseteq w_K$, so
$\sqrt{w_J} = \sqrt{w_K}.$
\end{proof}

It follows that if $H$ is itself pro$-p$, then $w_H \subseteq P_H$. If in addition $\Omega_{G/H}$ is semiprime, then $w_H = P_H$.

Our main result on localisable semiprime ideals is

\begin{thmd}
\label{PHLocalisable}
$P_H$ is classically localisable in $\Omega_G$.
\end{thmd}

In particular, if $G$ is itself a pro$-p$ group and $H$ is a
closed normal subgroup of $G$ with $\Delta^+ (G/H) = 1$, then $w_H$
is a prime localisable ideal of the local $\Fp -$algebra
$\Omega_G$ by Theorems A and D. One might imagine that in fact this is a special case
of a general phenomenon, namely that \emph{every} prime ideal $P$
of a local Noetherian $k$-algebra $R$ with $R/J(R) \cong k$ is
classically localisable. But this is not so - for a counterexample
see \cite[Example 7.2]{BHM}, an example first studied by Ramras.
However, in every such example known to us, there is a classically
localisable semiprime ideal $Q$ contained in $P$, with $P$ a
minimal prime ideal over $Q$ - that is, in the jargon of
\cite{Jat}, $P$ belongs to a finite clique. Indeed, the same is
true in all known examples when $R$ is assumed to be merely
semilocal rather than local.

\subsection{Finite subgroups} \label{finite}
We start the proof of Theorem D by taking care of a special case.

\begin{prop}
\label{SpecialCase}
Let $H$ be a \emph{finite} normal subgroup of the compact $p$-adic analytic group $G$. Then $P_H$ is classically localisable in $\Omega_G$.
\end{prop}
\begin{proof}
By Lemma \ref{PHWellDefined}, $P_H = \sqrt{w_{\{1\}}} = \sqrt{0}$
is the prime radical of $T=\Omega_G$. Pick an open normal uniform
subgroup $N$ of $G$. As in the proof of Proposition \ref{OuterActionSpecialCase}, let $R =
\Omega_N$ and $S = R \backslash \{0\}$; then $S$ is an Ore set in $T$ and
\[TS^{-1} \cong (RS^{-1})*(G/N)\]
is a crossed product of a division ring with a finite group. Thus
$TS^{-1}$ is Artinian, so any regular element in $TS^{-1}$ is a
unit by \cite[Proposition 3.1.1]{MCR}.

The natural map $T \to TS^{-1}$ is an injection, $TS^{-1}$ is a
partial ring of quotients of $T$, and every regular element of $T$
is a unit in $TS^{-1}$. Hence $TS^{-1}$ is an Artinian (right and
left) quotient ring of $T$, in the sense of 2.1.3 and 2.1.14 of
\cite{MCR}. Thus, the set $\mathcal{C}_T(0)$ of regular elements
of $T$ must be an Ore set in $T$. By Small's Theorem
\cite[Corollary 4.1.4]{MCR}, $\mathcal{C}_T(0) =
\mathcal{C}_T(\sqrt{0})$, so $\sqrt{0}$ is localisable in $T$ as
required. Since $\sqrt{0}TS^{-1}$ is the nilpotent radical of
$TS^{-1}$ it clearly has the Artin-Rees property, so $P_H$ is
classically localisable.
\end{proof}

\subsection{The Artin-Rees property} \label{AR}
The following result is inspired by Theorem 4.6 of \cite{Ven}.

\begin{lem}
\label{GradedRing} Let $N$ be a closed normal uniform subgroup of
the compact $p$-adic analytic group $G$ and let $E = E_G(N)$. Then the
graded ring of $\Omega_G$ with respect to the $w_N$-adic
filtration is isomorphic to
\[\Omega_{E/N}[X_1,\ldots,X_d] * (G/E) \]
where $d = \dim(N)$, and is a Noetherian ring.
\end{lem}
\begin{proof}
Let $\mathfrak{m} = w_{N,N}$ be the maximal ideal of $\Omega_N$,
so that $w_N = \mathfrak{m}\Omega_G$. Since $w_N^n =
\mathfrak{m}^n\Omega_G = w_{N,E}^n\Omega_G$, we see that
\[\gr_{w_N}\Omega_G \cong (\gr_{w_{N,E}} \Omega_E) \otimes_{\Omega_E}\Omega_G \]
as right $\Omega_G$-modules. Since $\gr_{w_{N,E}} \Omega_E$ is a
subring and since $E$ has finite index in $G$, the ring
$\gr_{w_N}\Omega_G$ is isomorphic to a crossed product of
$\gr_{w_{N,E}} \Omega_E$ with $G/E$. We may hence assume that $E =
G$, so $G$ acts trivially on $N/N^p$.

Next, $gN^p \mapsto g - 1 + \mathfrak{m}^2$ gives an isomorphism
of the $G$-modules $N/N^p$ and $\mathfrak{m}/\mathfrak{m}^2$, so
$G$ acts trivially on $\gr_{\mathfrak{m}}\Omega_N$. Now,
\[\gr_{w_N}\Omega_G \cong (\gr_{\mathfrak{m}} \Omega_N) \otimes_{\Omega_N}\Omega_G \cong (\gr_{\mathfrak{m}}\Omega_N)\otimes_{\Fp}\Omega_{G/N} \cong \Omega_{G/N}[X_1,\ldots,X_d]\]
as $\Omega_G$-modules, because $\gr_{\mathfrak{m}}\Omega_N \cong
\Fp[X_1,\ldots,X_d]$ by Theorem 7.24 of \cite{DDMS}. Since $G$
acts trivially on $\gr_{\mathfrak{m}}\Omega_N$, this is a ring
isomorphism.
\end{proof}

\begin{prop}
\label{WJArtinRees} If $J$ is a closed normal uniform subgroup of
$G$, then $w_J\triangleleft \Omega_G$ has the right and left
Artin-Rees properties.
\end{prop}
\begin{proof}
Let $N$ be an open normal subgroup of $G$. Since $N \cap J$ is
open in $J$ it has $p$-power index in $J$. Thus $NJ/N$ is a normal
$p$-subgroup of $G/N$, so $w_J^k\subseteq w_{NJ}^k \subseteq
w_{N}$ for some $k\geq 0$, by Lemma \ref{FNPSubgroup}. It follows
that the $w_J$-adic topology is finer than the defining topology
on $\Omega_G$ (where the base of open neighbourhoods of 0 is
$\{w_N : N\triangleleft_o G\}$), so $\Omega_G$ is complete with
respect to the $w_J$-adic filtration. Also, the graded ring of $\Omega_G$ with respect to the $w_J$-adic filtration is Noetherian by the lemma. By \cite[Chapter II, Proposition 2.2.1]{LVO}, this filtration is Zariskian. Hence, by
 \cite[Chapter II, Theorem 2.2]{LVO}, it has the Artin-Rees
property, which is easily seen to imply that the ideal $w_J$ has
the Artin-Rees property.
\end{proof}

Since it is trivial to check that if an ideal $I$ of a Noetherian
ring $R$ has the right Artin-Rees property then so does
$\sqrt{I}$, the following corollary is an  immediate consequence
of the proposition and the fact \cite[Corollary 4.3]{DDMS} that we
can find an open uniform subgroup of $H$ which is normal in $G$.

\begin{cor} \label{ARCor}
$P_H$ has the Artin-Rees property whenever $H$ is a closed normal
subgroup of $G$.
\end{cor}

\subsection{Proof of Theorem D} \label{mainpf} Let $H$ be a closed
normal subgroup of the compact $p-$adic analytic group $G$. Given
Corollary \ref{ARCor}, the theorem will follow if we show that
$P_H$ is localisable. As before, let $J$ be an open uniform
subgroup of $H$ which is normal in $G$. By Lemma
\ref{PHWellDefined}, $P_J = P_H$, so we may assume that $H$ is
actually uniform; in particular, $w_H \subseteq P_H =: P$. Since
$P$ has the Artin-Rees property, \cite[Theorem 4.2.10]{MCR} shows that it's sufficient to prove
that $P/P^n$ is localisable in $\Omega_G/P^n$ for all $n\geq 0$.

Now, the sets of ideals $\{w_{H,H}^k : k\geq 0\}$ and
$\{w_{K,H}:K\triangleleft_o H\}$ are cofinal in $\Omega_H$ by
\cite[Lemma 7.1]{DDMS}, so for all $n\geq 0$ there exists
$K\triangleleft_oH$ such that $w_{K,H} \subseteq w_{H,H}^n$. By
passing to an open subgroup, we may assume $K$ is normal in $G$.
Therefore, $w_K \subseteq w_H^n \subseteq P^n$, so it's sufficient
to show that $P_{H/K} = P/w_K$ is localisable in $\Omega_G / w_K
\cong \Omega_{G/K}$. This follows directly from Proposition
\ref{SpecialCase}.

\subsection{Centrally generated ideals} \label{central} Both in the proof of Theorem E
and when we come to consider the passage from localisability in
$\Omega_G$ to localisability in $\Lambda_G$ in $\S$\ref{lift} we
need to lift Ore sets over a centrally generated ideal. In fact,
for essentially no extra effort, we can lift over a
\emph{polycentral ideal} $I$ of the Noetherian ring $R$; that is,
$I = \Sigma_{i = 1}^n x_i R,$ where $x_1$ is central in $R$, and,
for $2 \leq j \leq n$, $x_j$ is central \emph{modulo}
$\Sigma_{i=1}^{j-1} x_i R.$

\begin{prop}
 Let $I$ be a polycentral ideal of the Noetherian ring $R$, with $I \subseteq P$, $P$ a semiprime ideal of $R$.
 Use $\overline{-}$ to denote images in $\overline{R} := R/I.$ Then
 $P$ is classically (right) localisable in $R$ if and
 only if $\overline{P}$ is classically (right) localisable in $\overline{R}$.
\end{prop}

\begin{proof} One direction is trivial. So suppose for the
non-trivial direction that $\overline{P}$ is classically right
localisable in $\overline{R}.$ By Noetherian induction we may
assume that $I$ is generated by a single central element $x$ of
$R.$ We first show that $P/x^nR$ is localisable in $R/x^nR$ for
all $n\geq 1$. By induction, assume that $n > 1$ and that
$P/x^{n-1}R$ is localisable in $R/x^{n-1}R$.

 Let $r \in R$ and $c\in\mathcal{C}_R(P)$. We can find $r' \in R$ and $c'\in\mathcal{C}_R(P)$ such that $rc' - cr' = x^{n-1}u$ for some $u \in R$.
 We can also find $u' \in R$ and $c''    \in\mathcal{C}_R(P)$ such that  $uc'' - cu' = xv$ for some $v \in R$.
 Hence, $rc'c'' - cr'c'' = x^{n-1}uc'' = x^{n-1}(cu'+xv)$ so $r(c'c'') - c(r'c'' + x^{n-1}u') = x^nv$ as required.

 Again, let $r \in R$ and $c\in\mathcal{C}_R(P)$. By \cite[Proposition 4.2.6]{MCR}, $xR$ has the right Artin-Rees property, so there exists $n$ such that
 $(rR + cR)\cap (xR)^n \subseteq rxR + cxR$. By the above, we can find $r' \in R$ and $c'\in\mathcal{C}_R(P)$ such that $rc' - cr' \in x^nR$, so there exist $a,b \in R$ such that
 $rc' - cr' = rxa + cxb$. Hence $r(c' - xa) = c(r' + xb)$ and the result follows since $x \in P$.

 Finally, the fact that the ideal $PR_P$ of $R_P$ has the
 Artin-Rees property follows from the facts that
 \begin{itemize}
 \item the image of $x$ in $R_P$ is central;
 \item $\overline{P}\overline{R}_{\overline{P}}$ has the
 Artin-Rees property in $\overline{R}_{\overline{P}}$,
 \end{itemize}
combined with the following statement: If an ideal $A$ of a
Noetherian ring $S$ contains a polycentral ideal $B$, such that
$A/B$ has the right Artin-Rees property in $S/B$, then $A$ has the
right Artin-Rees property in $S$. This last statement can be
proved by imitating the argument of \cite[Theorem 4.2.7(i)]{MCR}.
\end{proof}

The above proposition is a slight generalisation of a result of
McConnell, \cite{McC}, who dealt with the case $I = P$. See also \cite[Theorem 2.2]{PFSmith}.

\subsection{Proof of Theorem E} \label{corloc}For this we need a group-theoretic
result. Recall that a finite group is $p-$\emph{nilpotent} if it
has a normal subgroup of order prime to $p$, the factor by which
is a $p-$group.

\begin{prop}
Let $G$ be a compact $p-$adic analytic group such that
$\Delta^+(G)$ is a $p-$group. Suppose that every factor of $G$ by
an open normal subgroup is $p-$nilpotent. Then $G$ is a pro$-p$
group.
\end{prop}
\begin{proof}
Let $N$ be an open normal uniform subgroup of $G$. Then $N^{p^\epsilon}$ is also a normal subgroup of $G$, so $G/N^{p^\epsilon}$ is $p-$nilpotent. Hence we can find a normal $p'-$subgroup $K/N^{p^\epsilon}$ of $G/N^{p^\epsilon}$ such that $G/K$ is a $p-$group. Since $N/N^{p^\epsilon}$ is a $p-$group, we see that $K \cap N \leq N^{p^\epsilon}$. It follows that $[K,N] \leq N^{p^\epsilon}$, so $K \leq E_G(N)$ which is a pro$-p$ group by Lemma \ref{subpop}. Hence $K = N^{p^\epsilon}$ and therefore $G$ is a pro$-p$ group by \cite[Proposition 1.11(ii)]{DDMS}.
\end{proof}

The following lemma is presumably well known, but we are unable to find a reference for it in the literature.
\begin{lem}
Let $H$ be a finite normal subgroup of the compact $p-$adic
analytic group $G$. Then $\Omega_G$ is a free $\Omega_H-$module.
\end{lem}
\begin{proof}
Choose an open uniform normal subgroup $N$ of $G$. Since $N$ is
torsion-free, $N \cap H = 1$ so $HN \cong H \times N$. Now
$\Omega_{HN} = \Omega_{H \times N} \cong \Omega_H \otimes_{\Fp}
\Omega_N$ because $H$ is finite, so $\Omega_{HN}$ is a free
$\Omega_{H}-$module. Since $\Omega_{G}$ is a free
$\Omega_{HN}-$module, our result follows.
\end{proof}

\emph{Proof of Theorem E:}
Let $F = \Delta^+_{p'}(H)$. Since $p \nmid |F|$, it's well known that $w_{F,F}$ is generated by the central idempotent
\[f = 1 - \frac{1}{|F|}\sum_{g \in F} g,\]
so $w_F = f \Omega_G$ and $\Omega_G/f\Omega_G \cong \Omega_{G/F}$. Since $F$ is normal in $G$, $f$ is central in $\Omega_G$.

$\textit{(iii)} \Rightarrow \textit{(ii)}:$
By Proposition \ref{central} we may assume that $F = 1$, so $H$ is pro$-p$. Since $\Omega_{G/H}$ is semiprime, $w_{H,G} = P_H$ and the result follows from Theorem $D$.

$\textit{(ii)} \Rightarrow \textit{(i)}:$ Trivial.

$\textit{(i)} \Rightarrow \textit{(iii)}:$ Suppose that $w_{H,G}$
is localisable. Without loss of generality we can assume that $F = 1$.

Let $V$ be any open normal subgroup of $H$. We claim that $H/V$ is $p-$nilpotent. Since quotients of finite $p-$nilpotent groups are $p-$nilpotent, we may assume that $V$ is normal in $G$ by passing to an open subgroup. Writing $\overline{-}$ for images in $G/V$, $w_{\overline{H}, \overline{G}}$ is localisable in $\Omega_{\overline{G}}$.

Since $\Omega_{\overline{G}}$ is a free $\Omega_{\overline{H}}-$module by the lemma, it's easy to deduce that $w_{\overline{H}}$ is a localisable  maximal ideal of $\Omega_{\overline{H}}$. That is, the augmentation ideal of the group algebra $\Fp\overline{H}$ is localisable. Therefore, by \cite[Theorem 11.2.15]{Pass2}, $\overline{H}$ is $p-$nilpotent as claimed.

In particular, $\Delta^+(H)$ is $p-$nilpotent. Since $F = 1$, $\Delta^+(H)$ is a $p-$group. Therefore the proposition shows that $H$ is pro$-p$, completing the proof of the Theorem.

\section{Lifting information to $\Lambda_G$}\label{lift}

\subsection{} \label{lambda} Let $G$ be a compact $p$-adic analytic group. As was noted
in (\ref{Iwachat}), $\Lambda_G$ is semilocal with $p$ a central
regular element in $J(\Lambda_G)$, such that
\[\Lambda_G / p\Lambda_G \cong \Omega_G.\]
These properties establish a strong connection between $\Lambda_G$ and $\Omega_G$, and we
use them to deduce corresponding results for $\Lambda_G$ that we established for $\Omega_G$
in the previous sections.

\subsection{Primeness, semiprimeness and zero divisors}
\label{lamprime}The first and last parts of the following result
are simply stated for completeness and are due to Neumann, \cite{Neumann}. The proof of (i) was noted in
(\ref{sprime}); (iii) can if desired be deduced from Theorem C by
an easy lifting argument. Similarly, it is easy to deduce (ii)
from Theorem A, noting that $\cap_{i \geq 0} p^i \Lambda_G = 0$.

\begin{thmf}
Let $G$ be a compact $p$-adic analytic group. Then
\begin{enumerate}[{(}i{)}]
\item $\Lambda_G$ is semiprime;
\item $\Lambda_G$ is prime if and only if $\Delta^+(G) = 1$;
\item $\Lambda_G$ is a domain if and only if $G$ is torsionfree.
\end{enumerate}
\end{thmf}

\subsection{Localisation in $\Lambda_G$.} \label{lamloc} As an
immediate consequence of Proposition \ref{central} and Theorem
D, we can state

\begin{thmg}
Let $G$ be a compact $p$-adic analytic group and let $H$ be a
closed normal subgroup. Let $I_H$ be the inverse image of the
ideal $P_H$ of $ \Omega_G$ in $\Lambda_G$. Then $I_H$ is
classically localisable in $\Lambda_G$.
\end{thmg}

In view of \cite[Proposition 2.6]{CFKSV}, this is a direct
generalisation of \cite[Theorem 2.4]{CFKSV} from the case when
$G/H \cong \Zp$ to the case of an arbitrary closed normal subgroup
$H$ of $G$.

\begin{thmh}
Let $G$ be a compact $p-$adic analytic group and let $H$ be a
closed normal subgroup of $G$ such that $\Omega_{G/H}$ is semiprime.
Let $v_{H,G}$ denote the inverse image of $w_{H,G}$ in $\Lambda_G$. Then the following are equivalent:
\begin{enumerate}[{(}i{)}]
\item $v_{H,G}$ is localisable;
\item $v_{H,G}$ is classically localisable;
\item $H / \Delta^+_{p'}(H)$ is pro$-p$.
\end{enumerate}
\end{thmh}
\begin{proof}
Let us write $(i)_{\Lambda}, (ii)_{\Lambda},
(i)_{\Omega},(ii)_{\Omega}$ for the relevant statements in
Theorems H and E respectively. Using Proposition \ref{central} and
Theorem E it is easy to obtain the chain of implications
$$(ii)_{\Lambda} \Rightarrow (i)_{\Lambda}\Rightarrow (i)_{\Omega}
\Leftrightarrow (iii) \Leftrightarrow (ii)_{\Omega}\Leftrightarrow
(ii)_{\Lambda}, $$ proving the result.
\end{proof}

\section{Krull and global dimensions}

\subsection{} \label{LocalCrossedProduct} Let $G$ be a compact $p-$adic analytic group, $H$ a closed normal subgroup. By Theorem D, the semiprime ideal $P_H = P_{H,G}$ is always localisable in $\Omega_G$. We will denote the localisation of $\Omega_G$ at $P_{H,G}$ by $\Omega_{G,H}$. $\Lambda_{G,H}$ is defined analogously.

Now suppose that $K$ is an open uniform subgroup of $G$. Let $T/(H\cap K) = \Delta^+(K/(H\cap K))$; since $K/(H\cap K)$ is powerful, $K/T$ is uniform by \cite[Theorem 4.20]{DDMS}. Since $H \cap K$ is open in both $T$ and $H$, $\Omega_{G,H} = \Omega_{G,T}$ by Lemma \ref{PHWellDefined}, so we may assume that $H=T$. It follows from \cite[Proposition 4.31(ii)]{DDMS} that $H$ is uniform.

That is, we will assume for the remainder of this section that $G$
is a compact $p-$adic analytic group containing an open normal
uniform subgroup $K$ and a closed normal subgroup $H$ such that $H
\subseteq K$, and both $H$ and $K/H$ are uniform. Note that these
conditions imply $[H,K] \subseteq H \cap [K,K] \subseteq
H^{p^\epsilon}$, whence $E_K(H) = K$.

\begin{lem}
$\Omega_{G,H}$ is a crossed product of $\Omega_{K,H}$ with the finite group $G/K$:
\[\Omega_{G,H} = \Omega_{K,H} * (G/K).\]
There is an analogous statement for $\Lambda_{G,H}.$
\end{lem}
\begin{proof}
We have observed earlier that $\Omega_G = \Omega_K * (G/K)$. Let $S = \mathcal{C}_{\Omega_K}(P_{H,K})$, an Ore set in $\Omega_K$ by Theorem D. Since $K$, $H$ and $K/H$ are uniform, $P_{H,K} = w_{H,K}$ and $S = \Omega_G \backslash w_{H,K}$ consists of regular elements of $\Omega_K$. Since $H$ and $K$ are normal in $G$, $S$ is $G$-invariant and is hence an Ore set in $\Omega_G$ consisting of regular elements by \cite[Lemma 37.7]{Passman}. Moreover,
\[\Omega_G S^{-1} = \Omega_K S^{-1} * (G/K) = \Omega_{K,H} * (G/K).\]
Now, as in the proof of Proposition \ref{SpecialCase}, $\Omega_{G/H}$ has an Artinian quotient ring so $\mathcal{C}_{\Omega_{G/H}}(0) = \mathcal{C}_{\Omega_{G/H}}(\sqrt{0})$ by \cite[Corollary 4.1.4]{MCR}. Hence $T := \mathcal{C}_{\Omega_G}(P_{H,G}) = \mathcal{C}_{\Omega_G}(w_{H,G})$.  It's easy to check that $S \subseteq T$ so we have a ring homomorphism
\[ \Omega_{K,H} * (G/K) = \Omega_G S^{-1} \longrightarrow \Omega_G T^{-1} = \Omega_{G,H}.\]
To show that this map is an isomorphism, it's sufficient to show that $T \subseteq \Omega_G \subseteq \Omega_G S^{-1}$ consists of units in $\Omega_G S^{-1}$.
To this end, observe that
\[\Omega_G S^{-1} / w_{H,G} S^{-1} \cong (\Omega_{K/H} * (G/K)) \bar{S}^{-1} \cong D_{K/H} * (G/K)\]
is an Artinian ring, where $D_{K/H}$ is the division ring of fractions of $\Omega_{K/H}$ and $\bar{S} = \Omega_{K/H} \backslash \{0\}$.
Using \cite[Proposition 2.1.16(iv)]{MCR} it's easy to show that any element $x \in \Omega_G$ which is regular modulo $w_{H,G}$ is regular modulo $w_{H,G}S^{-1}$ inside $\Omega_G S^{-1}$. Since any regular element in an Artinian ring is a unit \cite[Proposition 3.1.1]{MCR}, we can find $y \in \Omega_G S^{-1}$ such that $xy \equiv 1 \mod w_{H,G}S^{-1}$.
Now, $w_{H,K}S^{-1}$ is the Jacobson radical of $\Omega_{K,H}$ so $w_{H,K}S^{-1}$ is contained in the Jacobson
radical $J$ of $\Omega_G S^{-1} = \Omega_{K,H} * (G/K)$ by \cite[Theorem 4.2]{Passman}.
Hence $w_{H,G}S^{-1} \subseteq J$, so $xy \equiv 1 \mod J$. It follows that $x$ is a
unit in $\Omega_G S^{-1}$ as required.

The argument to handle $\Lambda_{G,H}$ is similar.
\end{proof}

\subsection{}\label{FaithfullyFlat}
The following result will be very useful in the computation of Krull and global dimensions of $\Omega_{G,H}$.

\begin{lem}
$\Omega_{G,H}$ is a faithfully flat $\Omega_H-$module.
\end{lem}
\begin{proof}
By \cite[Lemma 4.5]{Bru} $\Omega_G$ is a flat $\Omega_H-$module.
Because localisation is flat, it follows that $\Omega_{G,H}$ is a
flat $\Omega_H$-module. Since $\Omega_H$ has a unique maximal
right (and left) ideal $\mathfrak{m}=w_{H,H}$, to complete the
proof it's sufficient to prove that $\mathfrak{m}\Omega_{G,H} \neq
\Omega_{G,H}$ by \cite[Proposition 7.2.3]{MCR}. But
$\mathfrak{m}\Omega_{G,H}$ is contained in the Jacobson radical of
$\Omega_{G,H}$ and is hence a proper ideal of $\Omega_{G,H}$, as
required.
\end{proof}

\subsection{Krull dimension}

To obtain the upper bounds in Theorem I we require

\begin{prop}
\label{Zariskian}
Let $J = w_{H,K}\Omega_{K,H}$ denote the unique maximal ideal of $\Omega_{K,H}$. Then the $J-$adic filtration on $\Omega_{K,H}$ is Zariskian and the graded ring of $\Omega_{K,H}$ with respect to this filtration is isomorphic to
\[D_{K/H}[X_1,\ldots,X_d]\]
where $D_{K/H}$ denotes the division ring of fractions of $\Omega_{K/H}$ and $d = \dim H$.
\end{prop}
\begin{proof}
Write $I$ for $w_{H,K}$. We have observed in the proof of
Proposition \ref{WJArtinRees} that the $I-$adic filtration on
$\Omega_K$ is Zariskian. Hence the Rees ring
\[\Rees{\Omega_K} = \cdots \oplus I^2t^{-2} \oplus It^{-1} \oplus \Omega_K[t]\]
is Noetherian. Since $S = \Omega_K \backslash I$ is an Ore set in $\Omega_K$, it's easily checked that $S \subseteq \Rees{\Omega_K}$ is an Ore set in $\Rees{\Omega_K}$. Hence
\[\Rees{\Omega_K}S^{-1} = \cdots \oplus J^2t^{-2} \oplus Jt^{-1} \oplus \Omega_{K,H}[t] \cong \Rees{\Omega_{K,H}}\]
is Noetherian. Since $J$ is the Jacobson radical of $\Omega_{K,H}$, the $J-$adic filtration is Zariskian.

The second assertion follows easily from Lemma \ref{GradedRing}, since $E_K(H) = K$.
\end{proof}

\begin{proof}[Proof of Theorem I]
This is now just a matter of putting all the pieces together. By
Lemma \ref{FaithfullyFlat} and \cite[Lemma 6.5.3(i)]{MCR},
$\Kr(\Omega_H) \leq \Kr(\Omega_{G,H})$. Now $\Kr(\Omega_{G,H}) =
\Kr(\Omega_{K,H})$ by Lemma \ref{LocalCrossedProduct} and
\cite[Proposition 10.1.11(ii)]{MCR}.

Since $\Kr(D_{K/H}[X_1,\ldots,X_d]) = d$ by \cite[Proposition 6.5.4(i)]{MCR}, part (i) follows from the proposition and \cite[Chapter II, Corollary 3.1.3]{LVO}.

Now, $\Kr(\Lambda_{G,H}) \leq \Kr(\Lambda_G)$, and this is finite
by \cite[Theorem A]{Ard}; moreover, the element $p$ of $I_H
\subseteq J(\Lambda_{G,H})$ is a central regular element of
$\Lambda_{G,H}$. Since $\Lambda_{G,H} / p\Lambda_{G,H} \cong
\Omega_{G,H}$, the second part follows from \cite[Theorem
1.8]{Walker}. (See also \cite[Theorem C]{Ard}.)
\end{proof}

\subsection{Homological properties} \label{hom}

Before embarking on the proofs, we recall that a Noetherian ring
$R$ is \emph{Auslander-Gorenstein} if it has finite right and left
injective dimensions, $n$ say, and the \emph{Auslander property}
holds for all finitely generated right and left $R-$modules $M$:
that is, for all integers $i$, $0 \leq i \leq n,$ and for all $j <
i$ and all submodules $N$ of $\mathrm{Ext}^i_R(M,R),$
$\mathrm{Ext}^j_R(N,R)= 0.$ An \emph{Auslander-regular} ring is
one which is Auslander-Gorenstein and has finite global dimension.

The following lemma is surely well-known, but we have been unable
to locate it in the literature. Note that it is actually concerned
with \emph{Frobenius extensions}, as defined, for instance, in
\cite{BeFar}. The proof of the lemma, coupled with the simple
application of it outlined in the proof of Theorem J, shows that a
Frobenius extension of an Auslander-Gorenstein ring is
Auslander-Gorenstein.

\begin{lem} Let $T = R*G$ be a crossed product of the ring $R$ with the
finite group $G$ and let $M_T$ be a right $T-$module. Then
\[\mathrm{Ext}_T^i(M,T) \cong \mathrm{Ext}_R^i(M,R)\]
as left $R-$modules, for all $n\geq 0$.
\end{lem}
\begin{proof} As functors from the category of right $T-$modules to the
category of left $R-$modules, $\mathrm{Ext}_R^i(-,R)$ and
$\mathrm{Ext}_T^i(-,T)$ are derived from $\Hom_R(-,R_R)$ and from
$\Hom_T(-,T_T)$ respectively. It's therefore sufficient to prove
that
\[\Hom_T(M_T,T_T) \cong \Hom_R(M_R,R_R)\]
as left $R-$modules. Since $M_R \cong M_T \otimes_T T_R$, it's
enough to show that $T \cong \Hom_R(T_R, R_R)$ as right
$T-$modules.

Define $\alpha : \Hom_R(T_R,R_R) \to T$ by $\alpha(f) = \sum_{g\in
G} f(\bar{g})\bar{g}^{-1}$. An easy check shows that $\alpha$ is a
right $T-$module map. Since $\beta : T \to \Hom_R(T_R,R_R)$ given
by $\beta(\sum_{g\in G} r_g \bar{g})(\bar{h}) = r_{h^{-1}}$ is an
inverse, the result follows.
\end{proof}

\begin{cor} Let $R = T*G$ be the crossed product of a ring $R$ with a finite group $G$. Then the injective dimensions of $R$ and $T$ viewed as right modules over $R$ and $T$ respectively, are equal.
\end{cor}

\begin{proof}[Proof of Theorem J]
We prove first the claims concerning global dimension.

(iii) Choose an open normal uniform subgroup $K$ of $G$ and continue with the reductions made in \ref{LocalCrossedProduct}, so in particular $H \subseteq K \subseteq G.$

Suppose that $G$ contains no elements of order $p$. By
\cite[Theorem 4.1]{Bru} and \cite[Corollaire 1]{Serre}, the global
dimension of $\Omega_G$ and hence of $\Omega_{G,H}$ is finite. By
Lemma \ref{LocalCrossedProduct}, $\Omega_{G,H}$ is the crossed
product of $\Omega_{K,H}$ with the finite group $G/K$, so
$\gld(\Omega_{G,H}) = \gld(\Omega_{K,H})$ by \cite[Lemma
2.2(i)]{Z}.

By Lemma \ref{FaithfullyFlat}, \cite[Theorem 7.2.6]{MCR} and
the results of Brumer and Serre \cite{Bru, Serre},
\[d:=\dim H = \gld(\Omega_H) \leq \gld(\Omega_{K,H}).\]

By Proposition \ref{Zariskian}, we know that $\Omega_{K,H}$ has a
Zariskian filtration with graded ring $D_{K/H}[X_1,\ldots,X_d]$,
so we may deduce from \cite[Theorem 7.5.3(iii)]{MCR} and
\cite[Chapter II, Theorem 3.1.4]{LVO} that
\[\gld(\Omega_{K,H}) \leq \gld(D_{K/H}[X_1,\ldots,X_d]) = d\]
and the first set of displayed values in (iii) follows. The second
display follows easily from \cite[Theorem 7.3.7]{MCR}, since
$\Lambda_{G,H} / p\Lambda_{G,H} \cong \Omega_{G,H}$.

(iv) Suppose that $1 \neq x \in G$ with $x^p = 1$ and $xH \in
\Delta^+ (G/H)$, and suppose for a contradiction that $\gld \Lambda_{G,H} < \infty$. Let $C = \{g \in G : [g,x] \in H\}$, the centralizer in $G$ of $xH$; this is an open subgroup of $G$. Choose a uniform subgroup $K$ of $C$ which is open and normal in $G$. By changing $H$ if necessary as in \ref{LocalCrossedProduct}, we can assume that $H \leq K$ and that both $H$ and $K/H$ are uniform. Moreover, we still have the property that $[K,x] \subseteq H$.

By Lemma \ref{LocalCrossedProduct}, we have
\[ \Lambda_{G,H} \quad \cong \quad \Lambda_{K,H}\ast (G/K). \]
Set $E = \langle K,x \rangle = K\rtimes \langle x \rangle$, so $\Lambda_{E,H} = \Lambda_{K,H}\ast\langle x \rangle$ is a skew group ring. By
\cite[Lemma 2.1(ii)]{Z} and the above isomorphism, $\gld
\Lambda_{K,H}\ast\langle x \rangle < \infty.$

Write $J$ for the Jacobson radical of $\Lambda_{K,H}$. Since $K/H$ is uniform, $\Lambda_{K,H}/J$ is isomorphic to the division ring of fractions $F$ of $\Omega_{K/H}$. Thus \cite[Theorem
4.5(ii)]{Z} implies that $F\ast \langle x \rangle$ is semisimple
Artinian. But $[K,x] \subseteq H$ and $F$ has characteristic $p$, so $x - 1$ is a non-zero
central nilpotent element of $F\ast \langle x \rangle$. This
contradiction shows that $\gld \Lambda_{G,H}$ must be infinite.

Finally, if $\gld \Omega_{G,H} < \infty$ then $\gld \Lambda_{G,H}
< \infty$ by \cite[Theorem 7.3.7]{MCR}, so this case follows as
well.

(i) Lemma \ref{LocalCrossedProduct} and the lemma proved above
show that it's enough to prove the result for $\Omega_{K,H}.$ For
the latter algebra, the Auslander-Gorenstein property is a
consequence of Proposition \ref{Zariskian} and \cite[Chapter III, Theorem 2.2.5]{LVO}, since a polynomial ring over a division ring is
Auslander regular by \cite[Chapter III, Theorem 2.3.5]{LVO}. The value of the injective dimension
follows from the global dimension calculation and the corollary above.

(ii) Filter $\Lambda_{G,H}$ using the powers of the ideal
$p\Lambda_{G,H}$. The associated graded algebra is isomorphic to
$\Omega_{G,H}[X]$, and hence is Auslander-Gorenstein by (i) and \cite[Chapter III, Theorem 2.3.5]{LVO}. Since the filtration is separated and complete, we deduce
that $\Lambda_{G,H}$ is also Auslander-Gorenstein using
\cite[Chapter III, Theorem 2.2.5]{LVO}. The value of the injective
dimension follows as in (i).
\end{proof}

Whilst the condition that $G$ has no elements of order $p$ is
certainly sufficient in order for $\Omega_{G,H}$ and
$\Lambda_{G,H}$ to have finite global dimension, it is not
necessary. To see this, take any compact $p-$adic analytic group
$G$ with elements of order $p$ but such that $\Delta^+(G) = 1$
(the wreath product of $\Fp$ with $\Zp$ is an example). Then
$\Omega_G$ is prime by Theorem A. Letting $H=1$, we see that
$\Omega_{G,H}$ is the simple Artinian ring of quotients of
$\Omega_G$ and as such has global dimension $0$. We believe that
in fact the necessary condition of Theorem J(iv) is also
sufficient:

\begin{conj}
The global dimension of $\Omega_{G,H}$ (and of $\Lambda_{G,H}$) is infinite if and only if there exists an element $x \in G$ of order $p$ such that the centralizer of $xH$ in $G$ is open in $G$.
\end{conj}

\end{document}